\documentclass[12pt]{article}

\usepackage{amsfonts}
\usepackage{amsthm}
\usepackage{amsmath}
\usepackage{amsfonts}
\usepackage{latexsym}
\usepackage{amssymb}

\newtheorem{fed}{Definition}[section]
\newtheorem{teo}[fed]{Theorem}
\newtheorem{lem}[fed]{Lemma}
\newtheorem{cor}[fed]{Corollary}
\newtheorem{pro}[fed]{Proposition}
\newtheorem{rem}[fed]{Remark}

\def\dem{{\it Proof.\ }\rm}

\newfont{\bb}{msbm10}

\def\hrarr^#1_#2{ \mathrel{
\mathop{\hbox to .3in{\rightarrowfill}}
\limits^{\scriptstyle#1}_{\scriptstyle#2}  }}

\def\hlarr^#1_#2{ \mathrel{
\mathop{\hbox to .5in{\leftarrowfill}}
\limits^{\scriptstyle#1}_{\scriptstyle#2}  }}

\def\({\left(}
\def\){\right)}
\def\[{\left[}
\def\]{\right]}

\def\bull{\vrule height 1.0ex width .4ex depth -.1ex }

\def\inc{\subseteq}
\def\QED{\hfill \bull}

\def\*a{\#\sb a}

\def\H{{\cal H}}

\def\cH{{\cal H}}

\def\cP{{\cal P}}

\def\cQ{{\cal Q}}

\def\cW{{\cal W}}

\def\cM{{\cal M}}
\def\cN{{\cal N}}

\def\oo{\oplus}

\def\noi{\noindent}

\def\bm{\left(\begin{array}}
\def\em{\end{array}\right)}
\def\ben{\begin{enumerate}}
\def\een{\end{enumerate}}
\def\beq{\begin {equation}}
\def\eeq{\end {equation}}
\def\barr{\begin{array}}
\def\earr{\end{array}}

\def\H{{\cal H}}
\def\lh{{L(\H)}}
\def\lh+{{\lh^+}}

\def\cZ{{\cal Z}}
\def\cM{{\cal M}}
\def\cN{{\cal N}}

\def\cU{{\cal U}}

\def\({\left(}
\def\){\right)}
\def\[{\left[}
\def\]{\right]}

\def\bull{\vrule height 1.0ex width .4ex depth -.1ex }

\def\inc{\subseteq}
\def\QED{\hfill \bull}

\def\*a{\#\sb a}

\def\H{{\cal H}}

\def\cH{{\cal H}}

\def\cP{{\cal P}}

\def\cQ{{\cal Q}}

\def\cW{{\cal W}}

\def\cM{{\cal M}}
\def\cN{{\cal N}}

\def\noi{\noindent}

\def\bm{\left(\begin{array}}
\def\em{\end{array}\right)}
\def\ben{\begin{enumerate}}
\def\een{\end{enumerate}}
\def\beq{\begin {equation}}
\def\eeq{\end {equation}}
\def\barr{\begin{array}}
\def\earr{\end{array}}

\def\H{{\cal H}}
\def\lh{{L(\H)}}
\def\lh+{{\lh^+}}

\def\cZ{{\cal Z}}
\def\cM{{\cal M}}
\def\cN{{\cal N}}

\def\cU{{\cal U}}

\def\cJ{{\cal J}}

\def\M{{\cal M}}
\def\N{{\cal N}}

\def\H{{\cal H}}
\def\XX{{\mathfrak{X}}}
\def\YY{{\mathfrak{Y}}}

\def\inc{\subseteq}
\def\QED{ \hfill \bull}

\def\H{{\cal H}}

\def\cH{{\cal H}}

\def\cP{{\cal P}}

\def\cQ{{\cal Q}}

\def\cM{{\cal M}}
\def\cN{{\cal N}}

\def\noi{\noindent}

\def\bm{\left(\begin{array}}
\def\em{\end{array}\right)}
\def\ben{\begin{enumerate}}
\def\een{\end{enumerate}}
\def\beq{\begin {equation}}
\def\eeq{\end {equation}}
\def\barr{\begin{array}}
\def\earr{\end{array}}

\def\H{\cH}

\def\cZ{{\cal Z}}
\def\D{{\cal D}}

\def\J{{\cal J}}

\def\cA{{\cal A}}
\def\cW{{\cal W}}

\def\XT{\mathfrak{X}_T}
\def\YS{\mathfrak{Y}_S}
\begin{document}

\title{ Products of orthogonal projections and polar decompositions}

\author{ G. Corach$^{*}$ and A. Maestripieri$^{*}$
\thanks{\noindent Partially supported by  PICT 5272 FONCYT and UBACYT I023.}} 
\date{}
\maketitle
\begin{small}
\noindent $^{*}$Departamento de Matem\'atica, Facultad de Ingenier\'\i a, UBA
 and Instituto Argentino de Matem\'atica - CONICET, Saavedra 15, Buenos Aires (1083), Argentina.
 
\noindent e-mail: gcorach@fi.uba.ar, amaestri@fi.uba.ar
\end{small}
\maketitle
\medskip
\vskip.5truecm

\medskip


\vspace{1cm}

\vskip.5truecm
\begin{abstract}
We characterize the sets $\XX$ of all products $PQ$, and $\YY$ of all products $PQP$, where $P,Q$ run over all orthogonal projections and we solve
the problems 
\newline $\arg\min\{\|P-Q\|: (P,Q) \in \cal Z\}$, for $\cal Z=\XX$ or $\YY.$ 
We also determine the polar decompositions and Moore-Penrose 
pseudoinverses of elements of $\XX.$ 

\end{abstract}

%XXXXXXXXXXXXXXXXXXXXXXXXXXXXXXXXXXXXXXXXXXXXXXXXXXXXXXXXXXXXXXXXXXXXXXXXXXXXXXXXXXXXXXXXXXX
\section{Introduction}
Let $\cH$ be a Hilbert space; denote by $L(\cH)$ the algebra of all bounded linear operators on $\cH$ and by $\cP$ the set of all orthogonal projections in $L(\cH)$: 
$\cP=\{P\in L(\cH): P^2=P=P^*\}.$ The main goal of this paper is the study of the sets
$$
\XX=\{PQ: P,Q \in \cP\}
$$
and
$$
\YY=\{PQP: P,Q \in \cP\}.
$$
\newline In general, an operator $T\in \XX$ admits many factorizations like $PQ$. Crimmins (see comments below) proved that if $T\in \XX$ then $T=P_{\overline{R(T)}}P_{{N(T)}^\perp}$ (hereafter, $P_\cM$ denotes the orthogonal projection onto the closed subspace $\cM$, and $R(B)$, $N(B)$ denote the range and nullspace of $B$, respectively, for every operator $B \in L(\cH)$).  We characterize the set $\XT=\{(P,Q):
P,Q\in \cP,\, T=PQ\}$ and prove that the distinguished pair $(P_{\overline{R(T)}},P_{N(T)^\perp}) \in \XT$ is optimal in several senses. We study a similar problem for each $S\in \YY$: we characterize the set $\YY_S=\{(P,Q): P,Q\in \cP, \, S=PQP \}$ and find all pairs $(P_0,Q_0)\in \YY_S$ such that $\|P_0-Q_0\|=$min$\{\|P-Q\|: (P,Q)\in \YY_S\}.$ We also study the polar decomposition of operators in $\XX$ and show that the Moore-Penrose pseudoinverse operation is a bijection between $\XX$ and the set $\tilde\cQ$ of all closed (unbounded) projections. This bijection explains the coincidence between the set of all partial isometries which appear in the polar decomposition of oblique (i.e., not necessarily orthogonal) projections and those which appear in the polar decomposition of operators of $\XX.$

Products of orthogonal projections have attracted the attention of mathematicians from many different areas as functional analysis, mathematical physics, signal processing, numerical analysis, statistics, and so on. We refer the reader to the classical papers by J. Dixmier \cite{[D1]}, \cite{[D2]}, S. N. Afriat \cite{[Af]}, C. Davis \cite{[Da]} and P. Halmos \cite{[H1]}, \cite{[H2]} and  recent surveys by A. Gal\' antai \cite{[G]} and A. B$\rm{\ddot{o}}$ttcher and I. M. Spitkovsky \cite{[BS]}, which contain a large bibliography and several historical remarks. To their list we add a few papers which are closer to our results. I. Vidav \cite{[V]} studied the polar factors of oblique projections, and obtained several results which we recently rediscovered in \cite {[CM]}. In a paper of H. Radjavi and J. P. Williams on products of selfadjoint operators \cite{[RW]} there is a proof of a theorem by T. Crimmins which characterizes  the operators of $\XX$ in the following concise way: if $T\in L(\cH)$ then $T$ belongs to $\XX$ if and only if $T^2=TT^*T$; Crimmins also exhibited, for such $T$'s, what we call the {\it canonical factorization} $T=P_{\overline{R(T)}}P_{{N(T)}^\perp}$.  In \cite{[S]} Z. Sebesty\'en found a condition on an operator $T$ defined on  a subspace of $\cH$ in order to be the restriction of an orthogonal projection. We prove here that Sebesty\'en's condition is equivalent to Crimmins'.  More recently, A. Arias and S. Gudder \cite{[AG]} studied, in the more general setting of von Neumann algebras, what they call {\it almost sharp effects}, and which are, precisely, operators like $PQP$, for $P,Q\in \cP$. These  effects play a role in some problems of quantum mechanics. They found a characterization of the set $\YY$, which is very useful in our approach. It should be mentioned that in a complete different setting, S. Nelson and M. Neumann \cite{[NN]} found, for matrices, a characterization of the spectrum of elements of $\XX.$ It turns out that their conditions can be easily translated to the Arias-Gudder's theorem. T. Oikhberg \cite{[O1]}, \cite{[O2]} proved many results on operators which can be factorized as finite products of orthogonal projections. We close these comments by mentioning that some modern approaches to Heisenberg uncertainty principle, like those of Donoho and Stark \cite{[DS]} and Havin and J$\rm{\ddot{o}}$ricke \cite{[HJ]} (see also the survey by Folland and Sitaram \cite{[FS]}) are based on the compactness and spectral properties of certain products $PQ$, where $P$ and $Q$ respectively  project onto time-limited and band-limited signals.

We describe the contents of the sections. 
Section 2 contains some preliminary results. In section 3 we study some properties of operators of $\XX$ and characterize the set $\XT$ for $T\in \XX$, and we prove that the canonical factorization $T=P_{\overline{R(T)}}P_{{N(T)}^\perp}$ is optimal in the following senses: if $T=P_\cM P_\cN$ for some closed subspaces $\cM$, $\cN$, then (1) $\overline{R(T)} \inc \cM$ and ${N(T)}^\perp \inc \cN$; (2) $\|(P_\cM - P_\cN)x\| \ge \|(P_{\overline{R(T)}}-P_{{N(T)}^\perp})x \|$ for all $x\in \cH$; and (3) if $R(T)$ is closed then $\|P_{\overline{R(T)}}-P_{{N(T)}^\perp} \| < \|P_\cM-P_\cN\|$ for every other $(P_{\cM},P_{\cN}) \in \XT.$ In section 4 we start the study of the set $\YY$, by solving the problem $\arg \min\{\|P-Q\|: (P,Q)\in \YY_S\}$ for each $S\in \YY$. We include a theorem, whose proof is due to T. Ando, which describes, for fixed $P,Q\in \cP$, the set $\{H\in \cP: (PHP)^2=PQP\}.$ Section 5 is devoted to polar decompositions of elements of $\XX$. We characterize the set $\J_{\XX}$ (resp., $\XX^+$) of isometric  (resp., positive) parts of operators in $\XX$. In particular, we prove that $\XX=\{V^2: V\in \J_{\XX}\}$, $\XX^+= \YY$ and the map $T \longrightarrow V$, where $V$ is the isometric part of $T$, is a bijection between $\XX$ and $\J_{\XX}.$ The situation for the positive parts is different: using the above mentioned theorem, we parametrize, for every $S\in \YY$, the set $\{T\in \XX: |T|=S\}.$ In the last section we prove that the Moore-Penrose pseudoinverse of $T\in \XX$ is a closed unbounded oblique projection, and conversely. Using some results of Ota \cite{[O]} on closed unbounded projections, we extend this well-known theorem of Penrose \cite{[Pen]} and Greville \cite{[Gr]}, who proved this result for matrices. 

As observed by the referees, the techniques of Dixmier, Afriat, Davis and Halmos, as recently surveyed by Gal\'antai \cite{[G]} and B$\rm{\ddot{o}}$ttcher and Spitkovsky \cite{[BS]}, can be used to prove most of our results. See also the paper by Amrein and Sinha \cite {[AS]}. We have chosen to use more elementary tools, but we collect in a final remark a description of them.

%XXXXXXXXXXXXXXXXXXXXXXXXXXXXXXXXXXXXXXXXXXXXXXXXXXXXXXXXXXXXXXXXXXXXXXXXXXXXXXXXXXXXXXXXXXXXXXXXXXXX
\section{Preliminaries}
The direct sum of two closed subspaces $\M$ and $\N$ of $\H$ such that $\M\cap\N=\{0\}$ is denoted  by $\M\dot +\N,$ 
and if $\M$ and $\N$ are orthogonal we write $\M\oo \N$. For $A\in L(\H)$, $P_A$  stands for the orthogonal projection onto $\overline{R(A)}$. 
Denote $Gr(\cH)$ the Grassmannian manifold of $\cH$, i.e., the set of all closed subspaces $\cM$ of $\cH.$

The {\it Friedrichs angle} between $\cM \in Gr(\cH)$ and $\cN \in Gr(\cH)$ is $\alpha(\cM,\cN)\in [0,\pi/2]$ whose cosine is
$$
c(\cM,\cN)=\sup \{|\langle m,n\rangle|:\,m\in\cM\ominus\cN,\,\|m\|\le 1,\,n\in \cN\ominus\cM,\,\|m\|\le 1\},
$$
where $\cM\ominus\cN=\cM\cap(\cM\cap\cN)^\perp$.

The  {\it Dixmier angle} between $\cM$ and $\cN$ is $\alpha_0(\cM,\cN)\in [0,\pi/2]$  whose cosine is
$$
c_0(\cM,\cN)=\sup \{|\langle m,n\rangle|:\,\,m\in\cM,\,\|m\|\le 1,\,n\in \cN,\,\|m\|\le 1\}.
$$
It is easy to see that $c_0(\M,\N)=\|P_\M P_\N\|$; we collect several well-known facts on $c$ and $c_0$. The proofs can be found in the survey by F. Deutsch \cite{[Deu]}.

\begin{teo}\label{m+n} Given $\M, \N \in Gr (\cH)$ the following statements hold:
\begin{enumerate}
\item $ c(\M,\N)<1$ if and only if $\M+\N$ is closed if and only if $R(P_\M(I-P_\N))$ is closed;
\item $c_0(\M,\N)<1 \iff \M\cap\N=\{0\}$ and $\M+\N$ is closed;
\item $c(\M,\N)=c(\M^\perp,\N^\perp),$ i.e.,  the Friedrichs angle between $\M$ and $\N$ coincides with that between $\M^\perp$ and $\N^\perp;$
in particular, $\M + \N$ is closed if and only if  $\M^\perp+\N^\perp$ is closed.
\end{enumerate}
\end{teo}

We will use  the well known 
{\it Krein-Krasnoselskii-Milman equality}  
\begin{equation}\label{krasno}
\|P-Q\|=\max\{\|P(I-Q)\|,\,\,\|Q(I-P)\|\},
\end{equation}
valid for all $P, Q \in \cP$ (see \cite{[KKM]}, \cite{[AkGl]}, \cite{[KYos]}). 

\begin{pro} \label{p-q} Given $P,Q \in \cP$, there are four possible cases for the norms involved in Krein-Krasnoselskii-Milman equality, namely:
\begin{enumerate}
\item $\|P-Q\|<1$ and, then, $\|P(I-Q)\|=\|Q(I-P)\|<1$;
\item $\|P-Q\|=\|P(I-Q)\|=1$ and $\|Q(I-P)\|<1$;
\item $\|P-Q\|=\|Q(I-P)\|=1$ and $\|P(I-Q)\|<1;$
\item $\|P-Q\|=\|Q(I-P)\|=\|P(I-Q)\|=1.$
\end{enumerate}
In terms of the ranges and nullspaces of $P,Q$, the four possibilities read as follows:
\begin{enumerate}
\item $R(P)\dot +N(Q)=N(P)\dot +R(Q)=\cH$ and the angles of both decompositions coincide;
\item $R(P)+N(Q)=\cH$, the sum is not direct and $N(P)+R(Q)$ is a proper closed subspace;
\item $N(P)+R(Q)=\cH$, the sum is not direct and $R(P)+N(Q)$ is a proper closed subspace;
\item $N(P)+R(Q)$ and $R(P)+N(Q)$ are proper subspaces of $\cH$.
\end{enumerate}\end{pro}

Recall the definition of the Moore-Penrose pseudoinverse $T^\dagger$ of $T\in L(\cH).$ This is an operator
with  domain $R(T)\oplus R(T)^\perp$ defined by $T^\dagger(Tx)=x$ if $x\in N(T)^\perp$ and $T^\dagger|_{R(T)^\perp} =0.$
The reader is referred to the original paper by Penrose \cite{[Pen]} or the book by Ben-Israel and Greville \cite{[BG]}
for properties and theorems on $T^\dagger.$ We will use without explicit mention that $T^\dagger$ is bounded if and only if $R(T)$ is closed.
Notice that $T^\dagger T$ and $TT^\dagger$ behaves in a different way: the first one is always bounded; indeed, it coincides with $P_{N(T)^\perp}$; however,
the second is defined, and behaves like a projection, on the domain of $T^\dagger.$

%XXXXXXXXXXXXXXXXXXXXXXXXXXXXXXXXXXXXXXXXXXXXXXXXXXXXXXXXXXXXXXXXXXXXXXXXXXXXXXXXXXXXXXXXXXXXXX

\section{The set of products $PQ$}

In this section we  study the sets
$$
\XX=\{PQ:\,\, P,\,Q\,\in \cP\},\,\,\,\,\XX_{cr}=\{ T\in\XX :\, R(T) \textrm{ is closed}\}.
$$

We start with a theorem that gives two alternative characterizations of  the elements of $\XX$. The first one is due to T. Crimmins ( item 2), see Radjavi and Williams \cite{[RW]}, Theorem 8. The second (item 3) is a rewriting of a result by Z. Sebesty\'en for suboperators, see \cite{[S]}, Theorem 1.

\begin{teo}For any $T\in L(\cH)$, the following assertions are equivalent:
\begin{enumerate}
\item  $T\in\XX$;
\item $T^2=TT^*T$;
\item $\|Tx\|^2=\langle Tx,x\rangle$, for all $x\in N(T)^\perp$.
\end{enumerate}

 In this case, $T=P_{\overline{R(T)}}P_{N(T)^\perp}=P_{\overline{R(T)}}P_{\overline{R(T^*)}}=P_{N(T^*)^\perp} P_{N(T)^\perp}.$ 
 \end{teo} 
We will refer to the factorization obtained in the above theorem as the  {\it canonical factorization} of $T$. 
\medskip

\dem  1 $\rightarrow$ 3: If $T\in\XX$ there exist $P,Q\in\cP$ such that $T=PQ$. Observe that $N(Q)\subseteq N(T)$ so that $N(T)^\perp \subseteq N(Q)^\perp$ and  then $QP_{N(T)^\perp}=P_{N(T)^\perp}$, or $Qx=x$, for all $x\in N(T)^\perp$. Therefore, if $x\in N(T)^\perp$, then
$\|Tx\|^2=\langle T^*Tx,\,x\rangle=\langle QPQ x,\,x\rangle= \langle PQ x,\,Qx\rangle=\langle T x,\,x\rangle$, as wanted.

3$\rightarrow$ 2: If $\|Tx\|^2=\langle Tx,x\rangle$, for all $x\in N(T)^\perp$, then $\langle Ty,Ty\rangle=\langle Ty,P_{N(T)^\perp}y\rangle$, for all $y\in\cH$, because $TP_{N(T)^\perp}=T$. Hence 
$\langle T^*Ty,y\rangle=\langle P_{N(T)^\perp}Ty,y\rangle$ for all $y\in\cH$, or 
$T^*T=P_{N(T)^\perp}T=T^\dagger T^2$. Therefore, multiplying by $T$ both sides of this equality, $TT^*T=TT^\dagger T^2$.
But observe that $TT^\dagger$ is the orthogonal projection onto $\overline {R(T)}$, restricted to $R(T)$, and $R(T^2) \subseteq R(T)$. Then  $TT^*T=T^2$.

2 $\rightarrow$ 1: If $TT^*T=T^2$ then multiplying by (the possibly unbounded operator) $T^\dagger$ both sides of this equality, 
we get $P_{N(T)^\perp}T^*T=P_{N(T)^\perp}T$, and taking adjoints $T^*TP_{N(T)^\perp}=T^*P_{N(T)^\perp}$. Multiplying by
${T^*}^\dagger$, we get $P_{N(T^*)^\perp}TP_{N(T)^\perp}=P_{N(T^*)^\perp}P_{N(T)^\perp}$. But using that $N(T^*)^\perp=\overline{R(T)}$ and that $T=P_{\overline{R(T)}}TP_{N(T)^\perp}$, it follows  the equality $T=P_{\overline{R(T)}}P_{N(T)^\perp}$ so that in particular $T\in\XX$.
\QED

\bigskip

It is obvious that $T^*\in\XX$ if $T\in\XX.$ By the formula $T=P_{\overline{R(T)}}P_{N(T)^\perp}$, it is clear that
$T$ is determined by the closed subspaces ${\overline{R(T)}}$ and $N(T)$.

\begin{teo}\label{crim} Every $T\in \XX$ has the following properties:
\begin{enumerate}
\item $\overline{R(T)}\cap N(T)=\{0\}$;
\item $\overline{R(T)}\dot + N(T)$ is dense;
\item $\overline{R(T)}\dot + N(T)=\cH$ if and only if $R(T)$ is closed.
\end{enumerate}

\end{teo}

\dem 1. Let $x\in \overline{R(T)}\cap N(T)$. Then $P_{N(T)^\perp}x=0$ and $x=P_{\overline{R(T)}}x$. Therefore, 
$0=P_{N(T)^\perp}x=P_{N(T)^\perp}P_{\overline{R(T)}}x=T^*x$ so that $x\in N(T^*)=R(T)^\perp$.  Thus, $x\in\overline{R(T)}\cap R(T)^\perp=\{0\}$.

2. If $T\in\XX$ then also $T^*\in\XX$. Applying  1 to $T^*$ we get $N(T^*)\cap\overline{R(T^*)}=\{0\}$, or
$R(T)^\perp\cap N(T)^\perp=\{0\}$. Taking orthogonal complements we get that $\overline{R(T)}\dot + N(T)$ is dense.

3. Recall from Theorem \ref{m+n}  that $\cM +\cN^\perp$ is closed if and only if $R(P_\cM P_\cN)$ is closed and apply this to
$\cM=\overline{R(T)}$, $\cN=N(T)^\perp$. Since $T=P_\M P_\N$, from 2 we get the result.\QED

\medskip

\begin{cor} \label {pq}For any $P$, $Q\in\cP$ there exists only two alternatives:
\begin{enumerate}
\item $R(PQ)$ is closed and $R(PQ)\dot + N(PQ)=\cH$; or
\item $R(PQ)$ is not closed and  $\overline{R(PQ)}\dot + N(PQ)$ is a proper dense subspace of $\cH$.
\end{enumerate}
\end{cor}

\medskip

The next result is a reformulation of the canonical factorization property.

\begin{teo} \label{XT} Let $T\in\XX$. There exists a factorization $T=P_\cM P_\cN$ such that $\cM\dot +\cN^\perp=\cH$ if and only if $R(T)$ is closed.  
In this case, there exists only one such factorization, namely $T=P_{R(T)} P_{N(T)^\perp}$, which corresponds to the 
decomposition $\cH=R(T)\dot+N(T)$.
\end{teo}
\dem Observe that, by Theorem \ref{crim}, if $R(T)$ is closed then ${R(T)}\dot + N(T)=\cH$ and  $T=P_{R(T)} P_{N(T)^\perp}$.

Conversely,  if $T=P_\cM P_\cN$ and $\cM\dot +\cN^\perp=\cH$, then in particular $\cM +\cN^\perp$ is closed and, therefore, 
$R(T)=R(P_\cM P_\cN)$ is closed (see \cite{[B2]} or \cite{[I]}).  The uniqueness follows from the general lemma below.\QED

\begin{lem} If $\cM\dot+ \cN=\cH$, $\cM_1\dot+ \cN_1=\cH$, $\cM\supseteq\cM_1$ and $\cN\supseteq\cN_1$ then $\cM=\cM_1$ and $\cN=\cN_1$.
\end{lem}
\dem Straightforward.\QED

\medskip

{\rem \rm  If $P,Q\in \cP$ and $R(PQ)$ is closed, Theorem \ref{XT} and Corollary \ref{pq} do not imply that  $R(P)\dot + N(Q)=\cH$; however, it does imply that the operator $T=PQ$ admits a factorization $T=P'Q'$ such that $R(P')\dot + N(Q')=\cH.$
} 

\bigskip
Our next result describes all factorizations $T=P_\cM P_\cN$ for a given $T\in\XX$ and shows that the canonical factorization is optimal,  
in the following 
two senses: (1) if $T=P_\cM P_\cN$  then $\cM\supseteq \overline{R(T)}$ and $\cN\supseteq  N(T)^\perp$ or equivalently
$P_\cM\ge P_{\overline{R(T)}}$ and $P_\cN \ge P_{N(T)^\perp}$; (2) if $T=P_\cM P_\cN$ 
then $\|(P_\cM -P_{\cN})x\|\ge \|(P_{\overline{R(T)}}-P_{N(T)^\perp})x\|$, for all $x\in\cH.$

\begin{teo} \label{PPQQ} Let $T\in\XX$ and $\cM$, $\cN\in Gr(\cH)$. Then $T=P_\cM P_\cN$ if and only if there exist $\cM_1$, $\cN_1\in Gr(\cH)$  such that
\begin{enumerate}
\item $\cM=\overline{R(T)}\oplus \cM_1$;
\item $\cN={N(T)^\perp}\oplus \cN_1$;
\item $\cM_1 \perp \cN_1$;
\item $\cM_1\oplus \cN_1\subseteq R(T)^\perp\cap N(T)$.
\end{enumerate}
\end{teo}
\dem By Crimmins' theorem, it holds $T=P_{\overline{R(T)}}P_{N(T)^\perp}$. If $T=P_\cM P_\cN$ then, in particular, $R(T)\subseteq \cM$ and, 
since $\cM$ is closed, $\overline{R(T)}\subseteq \cM$.
Analogously, $\cN^\perp=N(P_\cN)\subseteq N(T)$ and therefore $\cN\supseteq N(T)^\perp$. Thus, $\cM_1:= \cM\ominus\overline{R(T)}$ 
and $\cN_1:=\cN\ominus N(T)^\perp$ are well-defined and items 1 and 2 are verified.
Also, $\cM_1\subseteq R(T)^\perp$ and $\cN_1\subseteq N(T)$.

Now we compute $T=P_\cM P_\cN$, using the decompositions 1 and 2, and we get
$$
P_{\overline{R(T)}}P_{N(T)^\perp}=T=P_\cM P_\cN=(P_{\overline{R(T)}}+P_{\cM_1})(P_{N(T)^\perp}+P_{\cN_1})=
$$
$$
=P_{\overline{R(T)}}P_{N(T)^\perp}+P_{\overline{R(T)}}P_{\cN_1}+P_{\cM_1}P_{N(T)^\perp}+P_{\cM _1}P_{\cN_1}
$$
and, after cancellation,
\begin{equation}\label{M1}
P_{\overline{R(T)}}P_{\cN_1}+P_{\cM_1}P_{N(T)^\perp}+P_{\cM _1}P_{\cN_1}=0
\end{equation}

By multiplying at left equation (\ref{M1}) by $P_{\overline{R(T)}}$, we get $P_{\overline{R(T)}}P_{\cN_1}=0$, because
$\cM_1\perp \overline{R(T)}$. From here we deduce 
also that $\cN_1\subseteq R(T)^\perp$.

We have now 
\begin{equation}\label{M2}
P_{\cM_1}P_{N(T)^\perp}+P_{\cM _1}P_{\cN_1}=0
\end{equation}
and, by multiplying at right by $P_{N(T)^\perp}$ we get
\begin{equation}\label{M3}
P_{\cM_1}P_{N(T)^\perp}=0
\end{equation}
because
$\cN_1\perp N(T)^\perp$; thus, 
\begin{equation}\label{M4}
P_{\cM _1}P_{\cN_1}=0
\end{equation}
and also $\cM_1\subseteq N(T)$.
This completes the first part.

Conversely, if $\M_1$, $\cN_1$ satisfies 1-4 then
$$
P_\cM P_\cN=(P_{\overline{R(T)}}+P_{\cM_1})(P_{N(T)^\perp}+P_{\cN_1})=P_{\overline{R(T)}}P_{N(T)^\perp}=T,
$$
because all other products vanish.\QED

\medskip
\begin{cor}  Let $T\in \XX$. Then $T$ admits a unique factorization $T=P_\cM P_\cN$ if and only if $R(T)^\perp\cap N(T)=\{0\}.$
\end{cor}
\medskip

\begin{cor}Let $T\in\XX$. If $T=P_\cM P_\cN$ then $ \|(P_\cM -P_\cN)x\|\ge \|(P_{\overline{R(T)}}-P_{N(T)^\perp})x\|$
for all $x\in\cH$, that is $(P_\cM -P_\cN)^2\ge (P_{\overline{R(T)}}-P_{N(T)^\perp})^2$.
\end{cor}
\dem In fact, $P_\cM-P_\cN=(P_{\overline{R(T)}}-P_{N(T)^\perp})+(P_{\cM_1}-P_{\cN_1})$ and the images of both terms
are orthogonal so
$
\|P_\cM x-P_\cN x\|^2= \|P_{\overline{R(T)}}x-P_{N(T)^\perp}x\|^2+ \|P_{\cM_1}x-P_{\cN_1}x\|^2.
$\QED

\bigskip

In what follows, for each  $T\in\XX$ denote  $\XT:=\{(P,Q): T=PQ\} $.

\begin{teo}\label {min} Let $T\in\XX$. If $R(T)$ is not closed, then $\|P-Q\|=1$ for all $(P,Q)\in\XT$.
If $R(T)$ is closed, then $\|P_{R(T)}-P_{N(T)^\perp}\|<1$ and $\|P-Q\|=1$ for every other $(P,Q)\in\XT$.
\end{teo}
\dem If $R(T)$ is not closed, then by Theorem \ref{crim},  it follows that $\overline{R(T)}\dot+ N(T)$ is a dense proper subspace 
of $\cH$ and, therefore, by (\ref{krasno}) and Theorem \ref{m+n} $\|P_{\overline{R(T)}}-P_{N(T)^\perp}\|=1;$ by the corollary above it follows that $\|P-Q\|=1$ for all $(P,Q)\in\XT$.

If $R(T)$ is closed, then $\cH=R(T)\dot+N(T)$ then, by Theorem \ref{m+n}, $c(R(T), N(T))=c_0(R(T), N(T))=\|P_{R(T)}P_{N(T)}\|=\|P_{R(T)}(I-P_{N(T)^\perp})\|<1$. Also, $T^*$ has closed range and in the same way, we obtain  that $ \|P_{N(T)^\perp}P_{R(T)^\perp}\|<1$, but $ \|P_{N(T)^\perp}P_{R(T)^\perp}\|=\|(I-P_{R(T)})P_{N(T)^\perp}\|$. Applying (\ref{krasno}), we get $\|P_{R(T)}-P_{N(T)^\perp}\|<1$ .

Finally, according to Theorem \ref{XT}, it follows that $(P_{R(T)},P_{N(T)^\perp})$
is the only element of $\XT$ with that property. Thus, if $(P,Q)$ is another element of $\XT$ then $R(P)+N(Q)=\cH$ but the sum is
not direct. Therefore  $\|P-Q\|=1.$\QED

%XXXXXXXXXXXXXXXXXXXXXXXXXXXXXXXXXXXXXXXXXXXXXXXXXXXXXXXXXXXXXXXXXXXXXXXXXXXXXXXXXXXXXX
\section{The set of products $PQP$}

Denote $\YY=\{PQP:\,\,P,\,Q\in\cP\}$ and for $S\in\YY$ denote $\YY_S=\{(P,Q):\,\,S=PQP\}$. 
This section is devoted to the study of these sets, following the lines of the preceding section.
First, we describe the set $\YY_S$ for a given $S\in\YY$.

\begin{pro} \label{ys} The set $\YS$
is the disjoint union of all sets $\XT$, where $T\in \XX$ satisfies $TT^*=S$. 
\end{pro}
\dem If $(P,Q)\in \YY_S$, then $S=PQP$, $T:=PQ \in \XX$ and $(P,Q)\in \XT$. Conversely, if $(P,Q)\in\XT$ for some $T\in \XX$ such that $S=TT^*$, then $S=PQP$, i.e., $(P,Q)\in \YY_S$\QED

\bigskip
The set $\YY$ was completely described by Arias and Gudder \cite{[AG]}. They proved that a positive operator $A\in L(\cH)$ belongs to $\YY$ if and only
if  $A\le I$ and dim$\overline{R(A-A^2)} \le$ dim$N(A).$ (Indeed, they proved a more complete result, valid for von Neumann algebras; in the case of factors, their result has
the form we mentioned.)

Given $S\in\YY$, we compute the norm $\|P-Q\|$ for every $(P,Q)\in \YS.$
\begin{teo} \label{yys} Let $S\in \YY$. Then:
\begin{enumerate}
\item If $R(S)$ is not closed then  $\|P-Q\|=1$ for every pair  $(P, Q) \in \YY_S$.
\item If  $R(S)$ is closed,  then for each  pair $(P,Q)\in\YY_S$ and $T=PQ$ the following alternative holds: either $T=PQ$ is not the 
canonical factorization of $T$, and then 
$\|P-Q\|=1$, or $P=P_{R(T)}$ and $Q=P_{N(T)^\perp}$, in which case $\|P_{R(T)}-P_{N(T)^\perp}\|$ is a constant $<1$ which is 
independent of the factorization $S=TT^*$; more precisely, $\|P_{R(T)}-P_{N(T)^\perp}\|=\|P_{R(S)}-S\|^{1/2}$.
\end{enumerate}
\end{teo}
\dem Recall that for every operator $B\in L(\cH)$, it holds $R(B)$ is closed if and only if $R(BB^*)$ is closed if and only if $R(B^*B)$ is closed: in fact, by the  polar decomposition it follows  $R(B)=R((BB^*)^{1/2})$; therefore, $R(B)$ is closed if and only if $R((BB^*)^{1/2})$ is closed if and only if $R(BB^*)$ is closed. For $B^*B$ it suffices to replace $B$ by $B^*$, because $R(B)$ is closed if and only if $R(B^*)$ is closed.
Consider $S \in \YY$. 
1) If $R(S)$ is not closed then for every $T \in \XX $ such that $TT^*=S$, it holds that $R(T)$ is not closed; by Theorem \ref{min}, it follows
that $\|P-Q\|=1$ for every pair $(P,Q)\in \XT$ and so, by  Proposition \ref {ys}, the same is true for every $(P,Q)\in \YY_S$.

\medskip

2)If $R(S)$ is closed, fix $T \in \XX $ such that $TT^*=S$. By Theorem \ref{min}, $\|P-Q\|=1 $ for every pair $(P,Q)\in\XT$ except for  the 
canonical pair $(P_ {R(T)}, P_{N(T)^\perp})$, for which $\|P_{R(T)}-P_{N(T)^\perp}\| < 1$. 
Consider another $L \in \XX$ such that $LL^*=S$. We claim that $\|P_{R(T)}-P_{N(T)^\perp}\|=\|P_{R(L)}-P_{N(L)^\perp}\|<1$.  In order
to prove this assertion, we make a series of remarks.
\begin{enumerate}
\item Observe that $R(S)=R(T)=R(L)$; denote $P=P_{R(S)}$.
\item If $E,F\in\cP$  then from 1 of Propostion \ref{p-q}, 
%$\|E-F\| \le max \{\|E(I-F)\|, \|(I-E)F\|\}$; 
it easily follows that if $\|E-F\|<1$ then  
$\|E-F\|=\|E(I-F)\|=\|(I-E)F\|$.
\item Since $\|P - P_{N(T)^\perp}\|<1$, then $\|P - P_{N(T)^\perp}\|=\|P(I-P_{N(T)^\perp})\|=\|PP_{N(T)}\|$.
\item  Observe that  $S=TT^*=PP_{N(T)^\perp}P=P-PP_{N(T)}P$, so that $PP_{N(T)}P=P-S$. 
\end{enumerate}

Thus, by items (3) and (4),  it follows that
$\|P - P_{N(T)^\perp}\|^2=\|PP_{N(T)}\|^2=$ 

\noi $=\|PP_{N(T)}P\|=\|P-S\|.$
\QED

{\rem \rm \label{yys1} The proof above shows that, if $S\in \YY$ has a closed range, then the set $\YY_S$ is the union of two disjoint subsets,
say $\cU=\{(P,Q)\in \YY_S: R(P)\dot +N(Q)=\cH \}$ and $\cZ=\{(P,Q)\in \YY_S: R(P)+N(Q)=\cH \,\,{\rm and} \,\, R(P)\cap N(Q)\neq \{0\} \}.$ The functional 
$(P,Q) \to \|P-Q\|$ takes the constant values $\|P_{R(S)}-S\|^{1/2}$ on $\cU$ and  $1$ on $\cZ$, respectively.
}
\medskip

The following is a technical result which will be used later on:
\begin{lem} \label{PmenosA} Let $P\in\cP$ and  $0\le A\le P$, then the following identities hold:
$$
\overline{R(P-A)}=\overline{R(P-A^2)}=\overline{R(P-A^{1/2})}
$$
and
$$
\overline{R(A-A^2)}=\overline{R(A(P-A))}=\overline{R(P_A-A)}.
$$
\end{lem}

\dem  Observe that the operators $A$, $P-A$, $P-A^2$ and $P-A^{1/2}$ are positive and commute because of the monotonicity of the positive square root; and the same holds with $P_A$ instead of $P$.

Also, from $(P-A^2)=(P+A)(P-A)$ and $P+A$ invertible on $R(P)$ we get $N(P-A^2)=N(P-A)$. Taking the orthogonal complements we have $\overline{R(P-A^2)}=\overline{R(P-A)}$, and similarly $\overline{R(P-A)}=\overline{R(P-A^{1/2})}$. 

Observe that $PA=A=AP$ so $A-A^2=A(P-A)$. To prove that $\overline{R(A(P-A))}=\overline{R(P_A-A)}$, observe that $N(A(P-A))=N(A(P_A-A))=N(P_A-A)$ and take orthogonal complement.\QED

The next theorem gives the form of an orthogonal projection $Q$ in the presence of another orthogonal projection $P$, in terms of $2 \times 2$ matrices induced by the decomposition $R(P)\oplus N(P)= \cH$; this type of result appeared, in some form, in the above mentioned papers by Afriat, Davis, Halmos, Arias and Gudder,  Gal\' antai, and  B$\rm{\ddot{o}}$ttcher and Spitkovsky.

\begin{teo} \label{ando} \label{aando} Let $P$ and $Q$ be  orthogonal projections, then the matrix representation of $Q$, under the decomposition $R(P)\oplus N(P)= \cH$, 
is given by 
\begin{equation}\label{ando1}
Q = \bm {cc} A & A^{1/2}(P - A)^{1/2}U^*\\  UA^{1/2}(P - A)^{1/2} &
U(P - A)U^* + \hat{Q}\em,
\end{equation}
where $A=PQP$, $U$ is  a partial isometry  with initial space $\overline {R(A(P-A))}$ and final space $\cW\subseteq N(P)$ and  $\hat{Q}$ is an orthogonal projection with 
$R(\hat{Q}) \subset  N(P)\ominus R(U)$.

Conversely, given $P\in\cP$, $0 \leq A \leq P$ such that $\dim \overline{R(A(P-A))}\le \dim N(P)$, a partial isometry $U$ with initial space 
$\overline{R(A(P - A))}$  and final space $\cW\subseteq N(P)$ and an orthogonal projection 
$\hat{Q}$ with $R(\hat{Q}) \subseteq N(P)\ominus R(U)$ the right-hand side of (\ref{ando1}) gives an orthogonal projection.
\end{teo}

\dem Given $P,Q\in\cP$, consider the matrix representation of $Q$ in terms of $P$:
$$
 Q\ =\bm {cc} Q_{11}& Q_{12}\\ Q_{21}& Q_{22}\em 
$$

Write $A := Q_{11}$ and $B:=Q_{22}$. Since $Q\ge 0$, it follows that 
$$0\le A\le P, \,\,\,0\le B\le I-P\,\,\,\textrm{ and}\,\,Q_{12}^*=Q_{21}.$$

Since $Q^2 = Q$, we also have 
%%%%%(2)
\begin{equation}\label{doss}
Q_{12}Q_{21} = A(P - A) \quad {\rm and}\quad
AQ_{12} + Q_{12}B = Q_{12}
\end{equation}
Since $Q_{12}^*=Q_{21}$, from the first equality we get  
$$|Q_{21}|^2 = A(P - A)\,\,\,\textrm{ or }|Q_{21}| = A^{1/2}(P - A)^{1/2},$$
so, we can  conclude that there is an  isometry $U$ from 
$\overline {R(A^{1/2}(P-A)^{1/2})}=$

\noindent $\overline {R(A(P-A))}$ onto $\cW\subseteq N(P)$ such that 
$$
Q_{21} = UA^{1/2}(P- A)^{1/2} \quad {\rm and}\quad 
Q_{12} = A^{1/2}(P - A)^{1/2}U^*.
$$

But applying Lemma \ref{PmenosA},  $\overline {R(A(P-A))}=\overline {R(P_A-A)}$.

It follows from the second identity  of (\ref{doss}) that 
$$
AA^{1/2}(P - A)^{1/2}U^* + A^{1/2}(P - A)^{1/2}U^*B = A^{1/2}(P - A)^{1/2}U^*.
$$
Observe that $A^{1/2}(P - A)^{1/2}=A^{1/2}(P - A)^{1/2}P_A$, by Lemma \ref{PmenosA}; then
$$
0=A^{1/2}(P - A)^{1/2} [U^*B-(P_A-A)U^*]=A^{1/2}(P - A)^{1/2} [U^*B-(P_A-A)U^*];
$$
this implies $R(U^*B-(P_A-A)U^*)\subseteq \overline{N(A(P-A))}$. Since
$R(U^*B-(P_A-A)U^*)\subseteq \overline{R(A(P-A))}$, then we have 

$$
U^*B = (P_A - A)U^* \quad \hbox{and hence}\quad UU^*B = U(P_A - A)U^*=BUU^*.
$$
Since $UU^*=P_U$ is an  orthogonal projection and $P_UB = BP_U$, we get that   
$$
B = U(P_A - A)U^* + \hat{Q}
$$
where $\hat{Q}$ is an orthogonal projection with $R(\hat{Q}) \subseteq N(P)\ominus R(U)$. 
Observe that $UP=U(P_A+P_{R(P)\ominus\overline{ R(A)}})=UP_A$, because ${R(P)\ominus R(A)}\subseteq N(A)\subseteq N(P_A-A)=N(U)$.
Then $B = U(P - A)U^* + \hat{Q}$.
 Therefore we arrive at (\ref{ando1}).

It is immediate to see that for $0 \leq A \leq P$ satisfying the dimension condition, a partial isometry $U$ with 
initial space $\overline{R(A(P - A))}$ and final space $\cW\subset N(P)$ and 
an orthogonal projection $\hat{Q}$ with $R(\hat{Q})\subseteq N(P)\ominus R(U)$ the right-hand side of (\ref{ando1})
gives an orthogonal projection. This completes the proof.\QED

\medskip

 As a consequence we get the following dilation result (cf. Theorem 5 and Corollary 6 from \cite{[AG]}):
 \begin{cor}
 Given a positive contraction $A\in L(\cH)$, there exists $Q\in\cP$ such that  $A=P_AQP_A$ if and only if
 $\dim \overline{R(A-A^2))}\le \dim N(A)$.
  \end{cor}
\medskip
 The next result will be useful in a characterization of the set $\YY$ by means of the polar decomposition (see next section).
\begin{cor}  \label{corando} Given $P,Q\in\cP$, there exists $H\in\cP$ which is a solution of 
\begin{equation}\label{pqp1/2}
(PQP)^{1/2}=PXP.
\end{equation}
Moreover, all the orthogonal projections which are solutions of (\ref{pqp1/2}) are parametrized as
$$
H=\bm {cc} A & A^{1/2}(P - A)^{1/2}U^* \\
UA^{1/2}(P - A)^{1/2} & U(P - A)U^* + \hat{H}\em
$$
where $A=(PQP)^{1/2}$,
$U$ is a partial isometry 
with initial space $\overline{R(A(P - A))}$ and final space $\cW\subseteq N(P)$ and $\hat{H}$ is 
an orthogonal projection with $R(\hat{H}) \subseteq N(P)\ominus R(U)$.
\end{cor}
\dem Let $A=PQP$; by the proof of the above theorem, $\dim \overline{R(P_A-A)}\le \dim N(P)$.
Consider $A^{1/2}$, then $0\le A^{1/2}\le P$. 
Therefore, applying Lemma \ref{PmenosA}, $\dim \overline{R(P_{A^{1/2}}-A^{1/2})}=\dim \overline{R(P_A-A)}\le \dim N(P)$. Finally, applying Theorem \ref{aando},
the proof is complete.\QED

{ \rem \rm Observe that  the above theorem contains an alternative proof of the result by Arias and Gudder \cite{[AG]}
 mentioned before, in the setting of Hilbert spaces.
 
 In \cite{[NN]} Nelson and Neumann proved that a set $\{\lambda_1, ...\,,\lambda_n\}$ is the spectrum of a $n\times n$ matrix $B=PQ$, where $P,Q\in\cP$, if and only if $\sharp \{i\,:\, 0<\lambda_i<1\}\le \sharp \{i\,:\, \lambda_i=0\} $.
Since the spectrum of $PQ$ coincides with that of $PQP$ it follows that the result by Nelson and Neumann is the finite-dimensional version of the theorem of Arias and Gudder.}

%XXXXXXXXXXXXXXXXXXXXXXXXXXXXXXXXXXXXXXXXXXXXXXXXXXXXXXXXXXXXXXXXXXXXXXXXXXXXXXXX
\section{Polar decomposition of $PQ$}

The {\it polar decomposition} of an operator $C\in L(\cH)$ is a factorization $C=V_C|C|$, where $V_C$ is a partial isometry, $|C|=(C^*C)^{1/2}$ and $N(V_C)=N(C)$. It is well known that this factorization exists and is unique \cite{[RS]}.
Morever, $R(V_C)=\overline{R(C)}$, $V_C V_C^*=P_{\overline{R(C)}}$, $V_C^*V_C=P_{N(C)^\perp}$ and $C=|C^*|V_C$. In what follows, $V_C$ will be called the {\it isometric part} of $C$ and $|C|$ the {\it positive part} of $C$.

Given a subset $\cA$ of $L(\cH)$ we consider the set $ \cA^+$ (resp., $\cJ_\cA$) which consists of all positive
(resp., isometric)  parts of members of $\cA$.

In \cite{[CM]} we characterized $\cQ^+$, where $\cQ$ is the set of all idempotents in $L(\H)$ (notice that in \cite{[CM]}, we used  the more cumbersome notation $L(\cH)^+_\cQ$) and $\cJ_\cQ$. We apply now the results above and those of \cite{[CM]} to characterize $\XX^+$, $\XX_{cr}^+$, $\cJ_\XX$ and $\cJ_{\XX_{cr}}$.

In \cite{[CM]} there is a characterization of the set $\cJ_\cQ$ of all partial isometries of oblique projections. More precisely, it is proven that, for a given
$V\in \cJ$, there exists $E\in \cQ$ with polar decomposition $E=V|E|$ if and only if $VP_{R(V)}$ is a positive operator with range $R(V)$. In other terms, 
the restriction of $VP_{R(V)}$ to $R(V)$ is a positive invertible operator in $L(R(V)).$ The next result proves that the squares of such isometries
exhaustes the set $\XX_{cr}.$

\bigskip

\begin{teo} \label{XCR} $$\XX_{cr}=\{V^2:\,\, V\in \cal{J}_\cQ\}.$$
\end{teo}

\dem By \cite{[Gr]}, $T\in \XX_{cr}$ if and only if $T^\dagger\in\cQ$ so that we only need to prove that, if $E\in\cQ$ has polar decomposition $E=V|E|$ then 
$E^\dagger={V^*}^2$, and use the general fact that $V^*$ is the partial isometry of $E^*$ in its polar decomposition.
For $E^\dagger ={V^*}^2$, observe that $N(E)=N(V)$ and $R(E)=R(V)$ so that $E^\dagger=P_{N(E)^\perp}P_{R(E)}=P_{N(V)^\perp}P_{R(V)}=
(V^*V)(VV^*)$. By the characterization of $\cal{J}_\cQ$, it holds $VP_{R(V)}=(VP_{R(V)})^*=P_{R(V)} V^*$, so that
$V^2V^*=V{V^*}^2$. Then, $E^\dagger=V^*V{V^*}^2$. But, since $V^*$ is the Moore-Penrose pseudoinverse of $V$, it holds $V^*VV^*=V^*$. 
Thus, $E^\dagger={V^*}^2$. This proves the theorem. \QED

\bigskip
This result will be extended to the whole $\XX$, after the characterization of the set $\cJ_\XX$ in the next theorem.

\bigskip

Let $T\in\XX$ such that $T=PQ$ is the canonical factorization of $T$. Then the left polar decomposition of $T$ has the form 
\begin{equation}\label{polarXX}
T=(PQP)^{1/2}V_T.
\end{equation}

Now we characterize the set $ \cJ_\XX=\{V\in \cJ: {\rm there \,\,exists}\,\, T\in\XX \,\,{\rm such \,\,that} \,\,V=V_T\}$, i.e., the partial isometries of the
polar decompositions of elements of $\XX$.

\begin{teo} \label{JXX} Given  $V\in \cJ$, then $V\in \cJ_\XX$ if and only if $V^2V^*\ge 0$ and  $\overline{R(V^2V^*)}=R(V)$. In this case, it holds 
$\overline{R(V)\dot +N(V)}=\cH$
\end{teo}
\dem
Let $V\in \cJ_\XX$, then  there exists $T\in \XX$ such that $V=V_T.$ Let $T=PQ$ be the canonical factorization of $T.$ Recall that $P=P_{\overline{R(T)}}=P_{R(V)}$ and, by the definition of the polar decomposition, $R(V)=\overline{R(T)}$. Therefore, $V^2V^*=V(VV^*)=VP.$ But, from
(\ref{polarXX})
we get that ${(PQP)^{1/2}}^\dagger T=PV=V$ so that $V={(PQP)^{1/2}}^\dagger PQ$ and
then,  $VP={(PQP)^{1/2}}^\dagger PQP=(PQP)^{1/2}$.
Therefore
$$
VP=|\,T^*| \in L(\cH)^+.
$$
Moreover, $R(V^2V^*)=R(VP)=R(|\,T^*|)=R(T)$ so that $\overline{R(V^2V^*)}=R(V)$.

Conversely, suppose that $V\in \cJ$ satisfies that $V^2V^*=VP_{R(V)}\ge 0$ and that $\overline{R(VP_{R(V)})}=R(V)$. 
Let  $A=VP_{R(V)}$ and $T=P_{R(V)}P_{N(V)^\perp}\in\XX$. Since $A$ is positive, in particular $A=V^2V^*=V{V^*}^2$.
 Then $T=(VV^*)(V^*V)=V^2V^*V=VP_{R(V)}V(=V^2)=AV$ and this is the polar decomposition of $T$. In fact, observe that $TT^*=AVV^*A=AP_{R(V)}A=A^2$ so that
 $|\,T^*|=A$; also $V$ is  a partial isometry with final space $R(V)=\overline{R(V^2V^*)}=\overline{R(A)}=\overline{R(T)}$
and nullspace $N(V)=N(T)$: $N(V)\subseteq N(T)$ and if $Tx=0$ then $AVx=0$; therefore $Vx\in N(A)\cap \overline{R(V)}=N(A)\cap \overline{R(A)}=\{0\}$.

The last assertion, namely that $\cH=\overline{R(V)+N(V)}=\cH$ if $V\in \cJ_\XX$, follows directly from Theorem \ref{crim}, by observing that $R(V)=\overline{R(T)}$ and $N(V)=N(T).$
\QED

\medskip
Given $T\in\XX$ with polar decomposition $T=|T^*|V$ then $T=P_{R(V)}P_{N(V)^\perp}$ is the 
canonical factorization of $T$.  By the previous results, it also holds that $R(T)$ is closed if and only if $R(V)\dot +N(V)=\cH$.

We have proved that if $T=V^2$ for a given $V\in\cJ_\XX$, then $T\in\XX$ and $V$ is the partial isometry of $T$. Therefore:

\begin{cor} Consider the map $\alpha: \cJ_\XX\longrightarrow L(\cH)$, $\alpha(V)=V^2$. Then $\alpha$ is a bijection from 
$\cJ_\XX$ onto $\XX$. In particular, $\XX=\{V^2:\,\,V\in\cJ_\XX\}$.
\end{cor}
\dem If $V\in\cJ_\XX$ then, by Theorem \ref{JXX}, $V^2V^*\ge 0$; in particular, $V^2V^*=V{V^*}^2$. Then $T=(VV^*)(V^*V)\in\XX$; but $T=V{V^*}^2V=V^2V^*V=V^2$, so that $\alpha(V)=V^2\in\XX$. Let $T\in\XX$; if $V$ is the  isometric part of $T$ then, by Theorem \ref{JXX} again, we get $V^2V^*\ge 0$ and $T=P_{\overline{R(T)}}P_{N(T)}=(VV^*)(V^*V)=V{V^*}^2V=V^2V^*V=V^2=\alpha(V)$. Thus, the isometric part of $T$ is $V$, so that $\alpha$ is surjective and  $\alpha^{-1}(T)=V$.\QED

\medskip

The last Corollary extends our previous results Theorem \ref{XCR} and \cite{[CM]}, Theorem 5.2.

\medskip
\begin{teo} \label{piso} Let $V\in\cJ$. Then $V\in\cJ_\XX$ if and only if $V$ has a matrix representation, in terms of the decomposition
$\H=R(V)\oplus R(V)^\perp$,  of the type
\begin{equation}\label{uu}
V=\bm {cc} A & (P - A^2)^{1/2}U \\
0 & 0\em
\end{equation}
where $P=P_V$, $0\le A\le P$, $\overline{R(A)}=R(V)$, $\dim \overline{R(P-A^2)}\le \dim R(V)^\perp$ and 
$U$ is a partial isometry  with initial space contained in $R(V)^\perp$  and final space $\overline{R(P - A^2)}$.

\end{teo}

\dem If $V\in\cJ_\XX$ then there exists $T\in\XX$ such that $V=V_T$. In the same way as in Theorem \ref{JXX}, if $T=PQ$
is the canonical factorization of $T$ then
$$
VP=(PQP)^{1/2}=A,
$$
where $\overline{R(A)}=R(V)$ and, by  Theorem 5 and Corollary 6 of \cite{[AG]}, $A$ satisfies that $0 \leq A \leq P$ and  $\dim \overline{R(P-A)}\le \dim N(P)$. By Lemma \ref{PmenosA},  $\dim \overline{R(P-A^2)}\le \dim N(P)$. 

Therefore
$$
V=\bm {cc} A & V_{12} \\
0 & 0\em
$$
is the matrix of $V$. Since $VV^*=A^2+V_{12}V_{12}^*=P$, then 
$|\,V_{12}^*|=(P-A^2)^{1/2}$, so that $V_{12}= (P-A^2)^{1/2}U$, where $U$ is a partial isometry  with initial space contained in $R(V)^\perp=N(A)$  and final space $\overline{R(P - A^2)}$.

Conversely, if $V$ has the matrix representation (\ref{uu}), with $A$ and $U$ satisfying the hypothesis of the theorem, then $VV^*=A^2+P-A^2=P$, so that $V\in\cJ$, $VP=A\ge 0$, $\overline{R(A)}=R(V)$ by hypothesis. Therefore, applying Theorem \ref{JXX}, it follows that $V\in \cJ_\XX$.\QED

\medskip

We end this section with a characterization of  the set 
$$
 \XX^+=\{A\in L(\cH)^+: {\rm there \,\,exists}\,\, T\in\XX \,\,{\rm such \,\,that} \,\,A=|\,T^*|\},
 $$
  i.e., the positive parts of the
polar decompositions of elements of $\XX$.
\begin{pro}\label{XXYY}
$$
\XX^+=\YY.
$$
\end{pro}

\dem Let $A\in \XX^+$. Then there exists $T\in\XX$ such that $A=(TT^*)^{1/2}$. If $T=PQ$ is the canonical factorization of $T$, then $A=(PQP)^{1/2}$ and applying Corollary \ref{corando} there exists $H\in\cP$ such that $A=PHP$ so that $A\in \YY$.

Conversely, let $A\in\YY$. Then there exist $P,Q\in\cP$ such that $A=PQP$ and we can assume that $P=P_A$. By Theorem 5 and Corollary 6 of \cite{[AG]}, it follows that $0\le A\le P$, $\dim\overline{R(P-A)}\le \dim N(A)$ and, by Lemma \ref{PmenosA},
$\dim\overline{R(P-A)}=\dim\overline{R(P-A^2)}$. In this case $P$ and $A$ satisfy the conditions of  Theorem \ref{piso}
and we can construct an operator $T\in \XX$;  more precisely, consider $T=AV$ with 
$$
V=\bm {cc} A & (P - A^2)^{1/2}U \\
0 & 0\em
$$
where $U$ is a partial isometry  with initial space contained in $R(V)^\perp$  and final space $\overline{R(P - A^2)}$. Then $TT^*= A^2$ or $|\,T^*|=A$. Therefore $A\in\XX^+$.\QED

\begin{cor} Consider the map $\beta: \XX\longrightarrow\YY$, $\beta(T)=|\,T^*|$. Then the fibre of $A\in\YY$ is given by
$$
\beta^{-1}(\{A\})=\{T\in\XX:\,\,T=\bm {cc} A^2 & A(P - A^2)^{1/2}U \\
0 & 0\em\}
$$
where $P=P_{\overline{R(A)}}$, 
$U$ is a partial isometry  with initial space contained in $N(A)$  and final space $\overline{R(P - A^2)}$.
\end{cor}
\dem Apply Proposition \ref{XXYY}.\QED
%XXXXXXXXXXXXXXXXXXXXXXXXXXXXXXXXXXXXXXXXXXXXXXXXXXXXXXXXXXXXXXXXXXXXXXXXXXXXXXXXXXXXXXX
\section{On the Moore-Penrose pseudoinverse of $PQ$}

As mentioned in the Introduction, Penrose \cite{[Pen]} and Greville \cite{[Gr]} proved that the Moore-Penrose pseudoinverse of an idempotent matrix is a product of two orthogonal projections, and conversely. A proof of the next result, which extends their theorem to closed range operators in $\XX$, appears in \cite{[CM]}.

\begin{teo} \label{Greville} Let $T\in L(\cH)$. Then $T\in\XX_{cr}$ if and only if there exists $E\in\cQ$ such that $T=E^\dagger$. 
In symbols, $\XX_{cr}=\cQ^\dagger$.
\end{teo}

The generalization of Penrose-Greville theorem for operators  $T\in \XX$ with non-closed range forces the consideration of a certain class of  unbounded projections. We refer the reader to the paper \cite{[O]} for the properties of those projections which naturally appear in this context.
In what follows, we consider the set $\tilde \cQ$ of closed  unbounded projections, i.e.,  operators $E$ with a dense domain ${\cal D}(E)$ such that  ${\cal D}(E)=N(E)\dot  +R(E)$, $N(E)$ is closed,
$R(E)$ is closed in $\cH$ and $E(Ex)=Ex$ for all $x\in {\cal D}(E)$.

\begin{teo}If $T\in\XX$ then there exists a closed unbounded
projection $E:{\cal D}(E)\longrightarrow \cH$ such that $T=E^\dagger$. Conversely,  if $E$ is any closed unbounded projection then 
there exists an element $T\in\XX$ such that $E^\dagger=T$. Moreover, the map 
$T\longrightarrow T^\dagger$ from $\XX$ onto $\tilde \cQ$ is a bijection.
\end{teo}
\dem Suppose that $T\in\XX$.  Then (see, e.g.,  \cite{[BG]}) $E=T^\dagger$ is an unbounded pseudoinverse of $T$ with dense domain 
${\cal D}(E)=R(T)\oplus R(T)^\perp$, $R(E)=N(T)^\perp$ and $E$ verifies $TET=T$, in $\cH$,  and $ETE=E$ in ${\cal D}(E)$. Since 
$R(E)=N(T)^\perp$ we get 
\begin{equation}\label{*}
P_{N(T)^\perp}Ex=Ex,\,\,\,\forall x\in{\cal D}(E).
\end{equation}
It also holds that 
\begin{equation}\label{**}
EP_{\overline{R(T)}}x=Ex,\,\,\,\forall x\in{\cal D}(E).
\end{equation}

In fact, if $x\in{\cal D}(E)$ then $Ex=E(P_{\overline{R(T)}}x+P_{R(T)^\perp}x)=EP_{\overline{R(T)}}x$ because $P_{\overline{R(T)}}x\in R(T)$ and 
$R(T)^\perp=N(E)$.

Observe also that $R(E)=N(T)^\perp\subseteq {\cal D}(E)$:  if $x\in N(T)^\perp$ then 
$x=P_{\overline{R(T)}}x+P_{R(T)^\perp}x= P_{\overline{R(T)}}P_{N(T)^\perp}x+P_{R(T)^\perp}x= Tx+P_{R(T)^\perp}x$ so that $x\in{\cal D}(E)$. 
Therefore $E^2$ is well defined in ${\cal D}(E)$.

Finally,  for $x\in{\cal D}(E)$, we get 
$$
E^2x=EP_{N(T)^\perp}Ex= EP_{\overline{R(T)}}P_{N(T)^\perp}Ex= ETEx=Ex.
$$
Observe that the first equality follows from (\ref{*})  and the second from (\ref{**}), because $P_{N(T)^\perp}Ex\in {\cal D}(E)$. We have proved 
that $E^2=E$ in $ {\cal D}(E)$; $R(E)=N(T)^\perp$ and $N(E)=R(T)^\perp$, both closed subspaces. This proves that $E$ is an unbounded closed 
projection, see Lemma 3.5 of \cite{[O]}, namely $E=P_{N(T)^\perp //R(T)^\perp}$.

Conversely, suppose that $\cM$ and $\cN$ are closed subspaces such that $\cM\dot +\cN$ is a  dense subspace of $\cH$. Let 
$E:\cM\dot +\cN\longrightarrow \cM$ be  the (unbounded) projection with domain ${\cal D}(E)=\cM\dot +\cN$ onto $\cM$ with nullspace $\cN$.
We will show that the unbounded operator $E$ is the pseudoinverse of an element of $\XX$, namely, $E=(P_{\cN^\perp}P_\cM)^\dagger$: 
in fact, $P_\cM Ex=Ex$, for every $x\in{\cal D}(E)$ and $EP_\cM=P_\cM$, in $\cH$, because  $R(E)=\cM$. Also, 
$R(P_{\cN^\perp}P_\cM)=R(P_\cM-P_{\cN}P_\cM)\subseteq  \cM\dot +\cN\subseteq {\cal D}(E)$. Therefore $EP_{\cN^\perp}P_\cM$ is well defined for 
every $x\in\cH$ and $EP_{\cN^\perp}P_\cM=E(I-P_\cN) P_\cM=P_\cM$, then
\begin{equation}\label{ABA}
EP_{\cN^\perp}P_\cM=P_\cM.
\end{equation}
Consider $x\in R(P_{\cN^\perp}P_\cM)(\subseteq {\cal D}(E))$ then $x= P_{\cN^\perp}P_\cM y$, for $y\in\cH$. Using  equation (\ref{ABA}) we get
$P_{\cN^\perp}P_\cM E x=P_{\cN^\perp}P_\cM E (P_{\cN^\perp}P_\cM y)=P_{\cN^\perp}P_\cM y=x$, then
$$
P_{\cN^\perp}P_\cM Ex=x,
$$
for every $x\in R(P_{\cN^\perp}P_\cM)$.

On the other side, if $x\in R(P_{\cN^\perp}P_\cM )^\perp= N(P_\cM P_{\cN^\perp})=(\cN^\perp\cap\cM)\oplus  \cN\subseteq {\cal D}(E)$ then $x=y+z$, 
with $y\in \cN^\perp\cap\cM$ and $z\in\cN$,  so that $Ex=Ey=y$. Therefore,
$$ 
P_{\cN^\perp}P_\cM Ex=P_{\cN^\perp}Ex=P_{\cN^\perp}Ey=P_{\cN^\perp}y=0.
$$
This proves that 
\begin{equation}\label{BAB}
P_{\cN^\perp}P_\cM E=P_{\overline{R(P_{\cN^\perp}P_\cM)}},\,\,{\textrm {in }}{\cal D}(E).
\end{equation}

Equations (\ref{ABA}) and (\ref{BAB}) prove that $E^\dagger=P_{\cN^\perp}P_\cM \in \XX$.\QED
\medskip

\bigskip
{\rem \rm a) Observe that the domain ${\cal D}=R(T)\oplus R(T)^\perp$ of the operator $E=T^\dagger$ in  the above theorem can be also expressed as a (not necessarily orthogonal) direct sum of two closed subspaces, more precisely ${\cal D}=N(T)^\perp\dot +R(T)^\perp=R(E)\dot +N(E)$: we have already proved that $N(T)^\perp\subseteq {\cal D}$ so that  $N(T)^\perp\dot +R(T)^\perp\subseteq {\cal D}$; to prove the other inclusion we have to check that $R(T)\subseteq N(T)^\perp\dot + R(T)^\perp$: let $x\in R(T)$, then we can compute $T^\dagger x=Ex$ and $Ex\in N(T)^\perp$. Therefore $Ex=Ex+ (I-E)x \in N(T)^\perp + R(T)^\perp$.

b) Let $T\in\XX$ with polar decomposition $T=V|T|$. Let us consider the operator
with domain $\D=R(T)\oplus R(T)^\perp$, defined by
$$
E=|T|^\dagger V^*|_\D.
$$
Observe that $V:N(T)^\perp\longrightarrow \overline{R(T)}$ is unitary and, by construction of $V$, $V(R(|T|))=R(T)$.
Then,  $V^*(R(T))=R(|T|)$; also observe that $|T|^\dagger(R(|T|))= N(T)^\perp$. Therefore,
$E$ is well-defined and $E(\D)=N(T)^\perp$. 

If $x\in R(T)^\perp$ then $Ex=|T|^\dagger V^*x=0$ because
$R(T)^\perp=N(V^*)$. Let us see that $E$ is the identity on $N(T)^\perp$; we have to check that $N(T)^\perp\subseteq\D$:
if $x\in N(T)^\perp$ then $x= P_{R(T)}x+P_{R(T)^\perp}x=P_{R(T)}P_{N(T)^\perp}x+P_{R(T)^\perp}x=Tx+P_{R(T)^\perp}x\in \D$. Then
$$
Ex=|T|^\dagger V^*(Tx+P_{R(T)^\perp}x)=|T|^\dagger V^*Tx=|T|^\dagger V^*V|T|x=|T|^\dagger |T|x=P_{N(T)^\perp}x=x.
$$
Therefore, $E=P_{N(T)^\perp//R(T)^\perp}$, and its left "polar  decomposition" is 
$$
E=|T|^\dagger V^*|_\D.
$$

We can also consider $T=|T^*|V_T$ to obtain  the right "polar decomposition" of $E$ given by $E=V^*|T^*|^\dagger$, in
${\cal D}$.} 

c) Finally, observe that the Moore-Penrose pseudoinverses of positive parts of elements of $\XX$ are the positive parts of elements of $\tilde\cQ$, i.e. ${(\XX^+)}^\dagger=\tilde\cQ^+$. 

\medskip

In \cite{[CM]}  the set of isometric parts of bounded oblique projections is characterized. Using this characterization, together with  the construction of the (left) polar decomposition of elements of $\tilde \cQ$ as above, and the fact that if $T\in\XX$ then $T^*\in\XX$, we get the following result:

\begin{cor}
 $$
 \cJ_\XX=\cJ_{\tilde\cQ}
 $$
 and
 $$
 \cJ_{\XX_{cr}}=\cJ_{\cal Q}.$$
\end{cor}

One of the referees noticed that several results of this paper can be proven following the techniques used in the theorem known as Halmos' two projections theorem. We refer the reader to the paper \cite {[BS]} (Theorem 1.1) for a recent presentation of this theorem, with a historical notice about the mathematicians involved in the proof. Given $P=P_{\M}$ and $Q=P_{\N}$, decompose 
$$
\M=(\M\cap \N) \oo (\M\cap \N^\perp)\oo \M_0,
$$
and 
$$
\M^\perp=(\M^\perp\cap \N) \oo (\M^\perp\cap \N^\perp)\oo \M_1,
$$
for certain closed subspaces $\M_0\subseteq \M$ and $\M_1\subseteq \M^\perp$. Therefore,
$$
\cH=(\M\cap \N) \oo (\M\cap \N^\perp)\oo (\M^\perp\cap \N) \oo (\M^\perp\cap \N^\perp)\oo (\M_1\oo \M_0).
$$

The two projections theorem says that if $\M_0$ or $\M_1$ is non trivial then there exists a unitary operator $R: \M_1\longrightarrow \M_0$ and operators $S$ and $C$ acting on $\M_0$ such that $0\le S\le I$, $0\le C\le I$, $S^2+C^2=I$,
$N(S)=N(C)=\{0\}$ and 
$$
P=I_{(\M\cap \N)}\oo I_{(\M\cap \N^\perp)} \oo 0_{(\M^\perp\cap \N)}\oo 0_{(\M^\perp\cap \N^\perp)}\oo \bm {cc} I& 0\\ 0 & R^*\em \bm {cc} I& 0\\ 0 & 0\em \bm {cc} I& 0\\ 0 & R\em,
$$
$$
Q=I_{(\M\cap \N)}\oo 0_{(\M\cap \N^\perp)}\oo  I_{(\M^\perp \cap \N)}\oo 0_{(\M^\perp\cap \N^\perp)}\oo \bm {cc} I& 0\\ 0 & R^*\em \bm {cc} C^2& CS\\ CS & S^2\em \bm {cc} I& 0\\ 0 & R\em.
$$
As a consequence, 
$$
PQ= I_{(\M\cap \N)}\oo 0_{(\M\cap \N^\perp)} \oo 0_{(\M^\perp\cap \N)}\oo 0_{(\M^\perp\cap \N^\perp)}\oo \bm {cc} I& 0\\ 0 & R^*\em \bm {cc}C^2& CS\\ 0 & 0\em \bm {cc} I& 0\\ 0 & R\em,
$$
$$
PQP=I_{(\M\cap \N)}\oo 0_{(\M\cap \N^\perp)} \oo 0_{(\M^\perp\cap \N)}\oo 0_{(\M^\perp\cap \N^\perp)}\oo \bm {cc} I& 0\\ 0 & R^*\em \bm {cc}C^2& 0\\ 0 & 0\em \bm {cc} I& 0\\ 0 & R\em,
$$
and 
$$
P-Q=0\oo I\oo  -I\oo 0\oo \bm {cc} I& 0\\ 0 & R^*\em \bm {cc} I-C^2& -CS\\ -CS & -S^2\em \bm {cc} I& 0\\ 0 & R\em,
$$
with the obvious notation. Using these representations, one can find proofs of some of the theorems of the paper. We have chosen a different aproach which does not rely on the two projections theorem.

\medskip
{\bf Acknowledgement:} We are deeply grateful to Prof. T. Ando, for various discussions on central issues in the paper.
 We also thank the referees for many valuable comments, which improved the presentation of the original manuscript.
 We thank Prof. T. Oikhberg, who brought to our attention the paper by Arias and Gudder.
 
%XXXXXXXXXXXXXXXXXXXXXXXXXXXXXXXXXXXXXXXXXXXXXXXXXXXXXXXXXXXXXXXXXXXXXXXXXXXXXXXXXX
 
\end{document}